\documentclass[11pt]{amsart}

\usepackage{amssymb}

\evensidemargin\oddsidemargin

\begin{document}
\baselineskip=18pt
\setcounter{page}{1}

\newtheorem{Conju}{Conjecture 1\!\!}
\newtheorem{Conjd}{Conjecture 2\!\!}
\newtheorem{Theo}{Theorem\!\!}
\newtheorem{Rqs}{Remarks}
\newtheorem{Lemm}{Lemma\!\!}
\newtheorem{Def}{Definition\!\!}
\newtheorem{Prop}{Proposition\!\!}

\renewcommand{\theConju}{}
\renewcommand{\theConjd}{}
\renewcommand{\theDef}{}
\renewcommand{\theTheo}{}
\renewcommand{\theLemm}{}
\renewcommand{\theProp}{}

\def\a{\alpha}
\def\b{\beta}
\def\B{{\bf B}}
\def\C{{\mathcal{C}}}
\def\CC{{\mathbb{C}}}
\def\E{{\mathcal{E}}}
\def\Da{{\rm D}_\a}
\def\Dq{{\rm D}_q}
\def\EE{{\mathbb{E}}}
\def\elaw{\stackrel{d}{=}}
\def\eps{\varepsilon}
\def\fxy{f_{x,y}}
\def\G{\gamma}
\def\fa{f_\a}
\def\HH{{\mathbb{H}}}
\def\hS{{\hat S}}
\def\hT{{\hat T}}
\def\hX{{\hat X}}
\def\ii{{\rm i}}
\def\K{{\bf K}}
\def\L{{\mathcal{L}}}
\def\lbd{\lambda}
\def\lacc{\left\{}
\def\lcr{\left[}
\def\lpa{\left(}
\def\lva{\left|}
\def\NN{{\mathbb{N}}}
\def\pb{{\mathbb{P}}}
\def\R{{\mathcal{R}}}
\def\rl{{\mathbb{R}}}
\def\racc{\right\}}
\def\rcr{\right]}
\def\rpa{\right)}
\def\rva{\right|}
\def\tEm{{\tilde E_m}}
\def\tY{{\tilde Y}}

\newcommand{\fin}{\vspace{-0.4cm}
                  \begin{flushright}
                  \mbox{$\Box$}
                  \end{flushright}
                  \noindent}

\title[Diffusion hitting times and bell-shape]{Diffusion hitting times and the bell-shape}

\author[Wissem Jedidi and Thomas Simon]{Wissem Jedidi and Thomas Simon}

\address{Department of Statistic and Operation Research, College of Sciences, King Saud University, P.O. Box 2455, Riyadh 11451, Saudi Arabia. D\'epartement de Math\'ematiques, Universit\'e de Tunis El Manar, Facult\'e des Sciences de Tunis, 2092 - El Manar I, Tunis, Tunisia. {\em Email}: {\tt wissem\_jedidi@yahoo.fr}}

\address{Laboratoire Paul Painlev\'e, Universit\'e Lille 1, 59655 Villeneuve d'Ascq Cedex, France. Laboratoire de Physique Th\'eorique et Mod\`eles Statistiques, Universit\'e  Paris-Sud, 91405 Orsay Cedex, France. {\em Email}: {\tt simon@math.univ-lille1.fr}}



\begin{abstract} Consider a generalized diffusion on $\rl$ with speed measure $m$, in the natural scale. It is known that the conditional hitting times have a unimodal density function. We show that these hitting densities are bell-shaped if and only if $m$ has infinitely many points of increase between the starting point and the hit point. This result can be viewed as a visual corollary to Yamazato's general factorization for diffusion hitting times.
\end{abstract}

\maketitle

\section{Introduction and statement of the result}

This paper deals with a certain distributional property of hitting times for generalized diffusions. We use the standard notation for the latter processes, as described e.g.  in Chapter V of \cite{RW}. Let $m$ be a string, that is a right continuous non-decreasing function from $[-\infty, +\infty]$ to $[-\infty, +\infty],$ with $m(-\infty) =-\infty, m(+\infty) = +\infty$ and $m(0-) = 0.$ Set
$$r_1\, =\, \sup\{x<0, \; m(x) =-\infty\}, \qquad r_2\, =\, \inf\{x>0, \; m(x) =+\infty\},$$
and define the positive measure $m(dx)$ on $[-\infty, +\infty]$ by
$$m(dx) = dm(x)\;\;\mbox{on $(r_1, r_2)$}, \quad m([r_1, r_2]^c) = 0\quad \mbox{and} \quad m(\{r_1\}) = m(\{r_2\}) = +\infty.$$
Let $\{B_t, \, t\ge 0\}$ be a linear Brownian motion and $\{L_t^x, \, t\ge 0, \, x\in\rl\}$ be its local time. Introducing the additive functional
$$A_t \; =\; \int_\rl L_t^x\, m(dx)$$
and its right-continuous inverse $\tau_t = \inf\{ u > 0, \, A_u > t \},$ we define the process
$$X_t \; = \; B_{\tau_t}, \qquad 0\le t< \zeta,$$
with lifetime $\zeta = \inf\{ t > 0, X_t = r_1\;\mbox{or}\; r_2\}.$ Here and throughout, the notation $\inf\emptyset = +\infty$ is implicitly assumed. The process $X$ is strongly Markovian with state space $E_m = (r_1, r_2)\cap {\rm Supp}\, m,$ and we set $\pb_x$ for its law starting from $x\in E_m.$ The measure $m$ is called the speed measure of $X.$  As is well-known, the above time-changed Brownian motion allows to construct all linear diffusions on an interval up to some monotonous, scale transformation - see again Chapter V. 7 in \cite{RW}. One can also choose a speed measure with discrete support, which leads to a class of continuous time Markov chains called gap diffusions in the literature - see \cite{KW} for details and examples. The above generalized diffusion is chosen on the natural scale, which does not cause any loss of generality in our problem.

The hitting time of $y\in E_m$ is defined by
$$\tau_y \; =\; \inf\{ t > 0; \; X_t = y\}.$$
This definition extends to $y = r_i$ for $i =1,2$ if we suppose $\vert r_i\vert < \infty$ and $r_i\in {\bar E_m},$ setting $\tau_{r_i} = \lim_{y\to r_i} \tau_y.$ In this situation we adjoin $r_i$ to $E_m$ and denote by $\tEm$ the extended state space. Let $x\in E_m$ and $y\neq x\in\tEm$ be such that $\pb_x[\tau_y < \infty] > 0.$ Without loss of generality we will suppose $y > x.$ This paper deals with the conditional hitting time distribution
$$\pi_{xy} (dt)\; =\; \frac{\pb_x[\tau_y \in dt]}{\pb_x[\tau_y < \infty]}\cdot$$
It is well-known that the probability measure $\pi_{xy}$ is absolutely continuous, and we denote by $\fxy$ its density function on $(0, +\infty).$ It was shown in Theorem 1.2 of \cite{R} that this density is always unimodal, in other words that it has a unique local maximum.

In this paper, we are interested in the following refinement of unimodality for $\fxy$. A smooth density function defined on a real interval is said to be {\em bell-shaped} if all its derivatives vanish at both ends of the interval and if its $n-$th derivative vanishes exactly $n$ times, for all $n\ge 1.$ Setting $n =1$ shows that a bell-shaped density is strictly unimodal. For $n=2$ the bell-shape property entails, as for the familiar bell curve, that there is one inflection point on each side of the mode and that the second derivative is successively positive, negative, and positive. The visual meaning of the bell-shape for $n= 3$ or $4$ is less immediate and we refer to the introduction of \cite{TS1} for details and references. Some standard density functions like the Cauchy, the Gaussian, the Gumbel or the Student are bell-shaped, as can be seen by a direct analysis. Showing this property for less explicit densities is however a more demanding task.

Consider now the generalized inverse Gaussian distribution, with density
$$c\, x^{\lbd -1}\, e^{-(\chi x^{-1} + \psi x)}$$
over $(0,+\infty),$ where $c$ is the normalization constant and the variation domain of $(\lbd, \chi, \psi)$ is described in \cite{BBH} p. 49. If $\chi > 0,$ the above density is bell-shaped, as can again be seen by a direct analysis, and it is also a hitting density for some SDE driven by Brownian motion - see Theorem 2.1. in \cite{BBH}. If $\chi = 0,$ the above density is not bell-shaped, and it is the hitting density of a generalized diffusion only if $\lbd \le 1$ - see Corollary 2 in \cite{Y2}. In this case it is also completely monotone, hence the first passage density to the nearest neighbour of some continuous time Markov chain with discrete state space - see \cite{B} p. 147 and the references therein. If $\chi > 0$ the speed measure of the underlying diffusion has an everywhere positive density, whereas if $\chi = 0, \lbd \le 1,$ the involved speed measure is atomic. This example suggests a general relationship between the bell-shape property for the hitting times of a generalized diffusion and the support of its speed measure. This connection is illustrated by the following characterization.

\begin{Theo} The density $\fxy$ is bell-shaped iff $m$ has an infinite number of increase points on $(x,y).$
\end{Theo}

This simple criterion applies to all SDE's driven by Brownian motion, and also to more singular processes like those of Examples 2.12 and 2.13 in \cite{Fre}. There exists an immense literature on hitting times for real diffusions, however it seems that the above interesting distributional property has not been investigated as yet. Our method to prove the theorem relies on a general factorization by Yamazato \cite{Y2}, and a total positivity argument which the first author had used in a previous paper  \cite{TS1} on stable densities. In this respect, it is worth recalling that diffusion hitting times are always  infinitely divisible. The bell-shape property for all positive self-decomposable densities with infinite spectral function at zero is stated as an open problem in \cite{TS1} - see Conjecture 1 therein.

\section{Proof of the Theorem}

\subsection{Proof of the if part} Translating the measure $m$ if necessary, we may and will suppose $x = 0$ and $y > 0.$ We appeal to the factorization obtained in Theorem 1 of \cite{Y2}, which reads
$$\pi_{0y} \; =\; \mu_1\, \ast\, \mu_2,$$
where $\mu_1$ has Laplace transform
$$\int_0^\infty e^{-\lbd x}\, \mu_1 (dx)\; =\; \prod_i \lpa\frac{a_i}{a_i +\lbd}\rpa$$
for an at most countable family of positive parameters $\{ a_i\}$ which is either empty or increasing and such that $\sum a_i^{-1} < \infty$, and $\mu_2$ is an absolutely continuous probability with completely monotone density function. Observe that the probability distribution $\mu_1$ belongs to the so-called Bondesson class - see Chapter 9 in \cite{SSV}. The L\'evy measure of the infinitely divisible distribution $\mu_2$ also fulfils a certain condition which is described in the Theorem of \cite{Y3}, but this will not be used in the sequel. Let us finally refer to Section 15.2 in \cite{SSV} for a thorough presentation of this classical result.

From the proof of Theorem 1 in \cite{Y2} p. 155, we also know that
$$\int_0^\infty e^{-\lbd x}\, \mu_1 (dx)\; =\; \frac{y}{\psi (y, \lbd)}$$
where $\psi (z,\lbd)$ is the unique continuous solution to
\begin{equation}
\label{ODE}
\psi(z,\lbd)\; =\; z \; +\; \lbd\int_{[0,z)} (z-t)\psi(t,\lbd) dm (t), \qquad z >0.
\end{equation}
In particular, the family $\{-a_i\}$ is the set of zeroes of the function $\lbd \mapsto \psi(y, \lbd).$ We first show the intuitively obvious fact that under the assumption that $m$ has infinitely many points of increase on $(0,y),$ the family $\{a_i\}$ is infinite. This simple remark is more or less already included in \cite{SSV} p.204, but we write it down for completeness.

Fix $n\ge 1$ and let $0 < y_1 < \ldots < y_{n+1} = y$ be such that each interval $(y_{i}, y_{i+1})$ contains a point of increase of $m.$ For each $i\in [1,n]$ and $\lbd > 0$ we have
\begin{eqnarray*}
\psi(y_{i+1},\lbd) & = & \psi(y_i, \lbd)\; +\; \int_{y_i}^{y_{i+1}} \frac{\partial\psi}{\partial z} (z, \lbd)\, dz\\
& \ge & (y_{i+1} - y_i) \frac{\partial \psi}{\partial z}\, (y_i, \lbd)\\
& \ge & \lbd (y_{i+1} - y_i)\,\psi(y_i, \lbd)\,m(y_i, y_{i+1})
\end{eqnarray*}
where the inequalities follow from (\ref{ODE}) and the fact that $z\mapsto \psi(z, \lbd)$ is positive, increasing and convex on $(0,y).$ An induction argument shows that there exists a positive constant $K$ such that $\psi(y, \lbd)  \ge  K \lbd^n$ for all $\lbd > 0.$ This being true for every $n\ge 1,$ we see that $\lbd \mapsto \psi(y,\lbd)$ grows faster than any polynomial as $\lbd \to +\infty,$ and hence that the family $\{a_i\}$ must be infinite.

We deduce that $\pi_{0y}$ is the probability distribution of the independent sum
$$Y\;\elaw\;\sum_{i\ge 1}\; Z_i \; +\; Z,$$
where $Z_i \sim {\rm Exp} (a_i)$ for all $i \ge 1$ with $0 < a_1 < a_2 <\ldots < a_n < \ldots,$ and $Z$ has a completely monotone density. This entails in passing that $f_{0y}$ is real analytic on $(0, +\infty)$ - see also \cite{P} for a weaker result in the framework of stochastic differential equations. Comparing the above infinite sum with the sum of $n+1$ independent copies of $Z_{n+1}$ also shows that $\pi_{0y}(0,z) = O(z^n)$ for every $n\ge 1$ as $z\to 0+$ (this latter fact is known in the literature as Ray's estimate, see Problem 4.6.5. p.134 in \cite{IM}) and hence that $f_{0y}^{(n)} (0+) = 0$ for every $n\ge 0.$

The remainder of the proof is close to that of the Theorem in \cite{TS1}. Introduce the notation
$$\pm^n\; =\; \lacc \begin{array}{cl}+ & \mbox{if $n$ is even}\\
- & \mbox{if $n$ is odd}\end{array}\right.$$
for every $n\in\NN.$ For $\{\eps_1, \ldots, \eps_n\}$ some finite sequence in $\{-,0,+\},$ we say that a smooth function $f :\, (0,+\infty) \to\rl$ is of type $\eps_1\ldots\eps_n$ if it has limits (finite or infinite) at zero and at infinity, vanishes on a finite set, and if the ordered sequence of its signs on $[0,+\infty]$ is given by $\{\eps_1, \ldots, \eps_n\}.$ For brevity, we will write $f\sim\eps_1\ldots\eps_n$ to express this property. With this notation, we need to show that
\begin{equation}
\label{Main}
f_{0y}^{(n)}\;\sim\; 0\!\pm^0\!0\!\pm^1\!0\cdots\pm^n\!0, \qquad n\ge 1.
\end{equation}
Observe first that this property holds true for $n = 1,$ by R\"osler's result and the fact that $f_{0y}'$ has isolated zeroes on $(0, +\infty).$ To show the property for all $n\ge 2$ we will use the same argument as in \cite{TS1}. Fix $n\ge 2$ and consider the independent decomposition
$$Y\; \elaw\; Y_n\; +\; \tY_n,$$
with 
\vspace{-0.25cm}
$$Y_n \; \elaw \; \sum_{i= 1}^{n+2}\; Z_i \; +\; Z\qquad\mbox{and}\qquad \tY_n\; \elaw\; \sum_{i\ge n+3}\; Z_i.$$
Setting $g_n$ and $h_n$ for the respective smooth densities of $Y_n$ and $\tY_n,$ we borrow from the Proposition in \cite{TS1} the following key-result.
\begin{Prop}[\cite{TS1}] With the above notations, one has
\begin{equation}
\label{gn}
g_n^{(i)} \; \sim \;  0\!\pm^0\!0\!\pm^1\!0\cdots\pm^i\!0, \qquad 0\le i \le n+1.
\end{equation}
\end{Prop}
In particular, the function $g_n$ has a $\C^{n+1}$ extension on $\rl$ and we can differentiate the convolution product to obtain
$$f_{0y}^{(n+1)}(u)\; =\; \int_0^\infty g_n^{(n+1)}(v) h_n(u-v) dv, \qquad u >0.$$
Using now the standard notation of \cite{K}, set $S^{-} (f)$ resp. $S^{+}(f)$ for the number of sign changes of a real function $f$ on $(0, +\infty),$ the zero terms being discarded resp. included. By (\ref{gn}), Theorem 3.1.(a) p. 21 in \cite{K} and the fact that the function $h_n$ is a P\'olya frequency function of infinite order - see e.g. Example 3.2.2 in \cite{B}, we have
$$S^{-} (f_{0y}^{(n+1)})\; \le\; S^{-} (g_n^{(n+1)})\; =\; n+1.$$
This inequality entails first that $f_{0y}^{(n)}(+\infty) = 0.$ Indeed, it is easily seen from the fact that $f_{0y}$ is a density function that $0$ must be a limit point of $f_{0y}^{(n)}$ at infinity. If it is not a true limit, then Rolle's theorem entails that $S^{-} (f_{0y}^{(n+1)}) = +\infty,$ a contradiction.

Since $f_{0y}^{(n)}(0+) = f_{0y}^{(n)}(+\infty) = 0$ for all $n\ge 0,$ an induction based on Rolle's theorem shows first that $S^{+}(f^{(n)}_{0y})\ge n$ for all $n\ge 1.$ On the other hand, since $f_{0y}^{(n)}$ has isolated zeroes on $(0,+\infty)$ and vanishes at zero and infinity, again Rolle's theorem shows that
$$S^+(f_{0y}^{(n)})\; \le\; S^-(f_{0y}^{(n+1)})-1\; \le\; n$$
for all $n\ge 1.$ This shows $S^+(f_{0y}^{(n)}) = n$ for all $n\ge 1,$ which readily entails (\ref{Main}) as required.
\subsection{Proof of the only if part} Again, we can suppose $x =0,$ translating $m$ if necessary. If $m$ has exactly $n$ points of increase on $(0,y)$ then an induction based on (\ref{ODE}) shows that $\lbd \mapsto \psi(y,\lbd)$ is a polynomial of degree $n$ and hence, that $\pi_{0y}$ is the probability distribution of the finite independent sum
$$Y\;\elaw\;\sum_{i= 1}^n\; Z_i \; +\; Z,$$
with the above notation. It follows then from the Proposition in \cite{TS1} that
$$(-1)^{n+i} f_{0y}^{(i)}\;\sim\; \pm^0 0\!\pm^1\!0\cdots\pm^n\!0$$
for all $i\ge n,$ so that $f_{0y}$ is not bell-shaped.
\subsection{Final remark} When $m$ has exactly $n$ points of increase on $(x,y)$, the Proposition in \cite{TS1} also provides some visual insight on the density $\fxy$.  When $n=1$ for example, it shows that $\fxy$ is "whale-shaped" (concave and then convex, all derivatives vanishing once).



\begin{thebibliography}{10, french}

\bibitem{BBH}
O.~Barndorff-Nielsen, P.~Bl\ae sild and C. Halgreen. First hitting time models for the generalized inverse Gaussian distribution. {\em Stoch. Proc. Appl.} {\bf 7}, 49-54, 1978.

\bibitem{B}
L.~Bondesson. {\em Generalized Gamma convolutions and related classes of distributions and densities.} Lect. Notes Stat. {\bf 76}, Springer-Verlag, New York, 1992.

\bibitem{Fre}
D.~Freedman. {\em Brownian motion and diffusion.} Holden-Day, San-Francisco, 1971.

\bibitem{IM}
K.~It\^o and H.~P.~McKean. {\em Diffusion processes and their sample paths.} Springer-Verlag, Berlin, 1965.

\bibitem{K}
S.~Karlin. {\em Total positivity. Vol I.} Stanford University Press, Stanford, 1968.


\bibitem{KW}
S.~Kotani and S.~Watanabe. Krein's theory of strings and generalized diffusion processes. In: Fukushima (ed.) {\em Functional analysis in Markov processes}, Lect. Notes in Math {\bf 923}, 235-259, 1982.

\bibitem{P}
E.~J.~Pauwels. Smooth first-passage densities for one-dimensional diffusions. {\em J. Appl Probab.} {\bf 24} (2), 370-377, 1987.

\bibitem{RW}
L.~C.~G.~Rogers and D.~Williams. {\em Diffusions, Markov processes and martingales. Vol 2: It\^o calculus.} Cambridge University Press, Cambridge, 1987.

\bibitem{R}
U.~R\"osler. Unimodality of passage times for one-dimensional strong Markov processes. {\em Ann. Probab.} {\bf 8} (4), 853-859, 1980.

\bibitem{SSV}
R.~L.~Schilling, R.~Song and Z.~Vondra{\v{c}}ek. {\em Bernstein functions.} De Gruyter, Berlin, 2010.

\bibitem{TS1}
T.~Simon. Positive stable densities and the bell-shape. {\em Proc. Amer. Math. Soc.}
{\bf 143} (2), 885-895, 2015.

\bibitem{Y2}
M.~Yamazato. Hitting time distributions of single points for 1-dimensional generalized diffusion processes. {\em Nagoya Math. J.} {\bf 119}, 143-172, 1990.

\bibitem{Y3}
M.~Yamazato. Characterization of the class of hitting distributions of 1-dimensional generalized diffusion processes. In: Shiryaev et al. (eds.) {\em Probability Theory and Mathematical Statistics. Proceedings of the 6th USSR-Japan symposium, Kiev, 1991}, 422-428, 1992.

\end{thebibliography}
\end{document}